\documentclass[11pt]{amsart}
\usepackage{amssymb,amsfonts,amscd,amstext}
\usepackage{epsfig,psfrag,graphicx}
\usepackage{graphicx,psfrag,subfig}
\newtheorem{theorem}{Theorem}[section]
\newtheorem{lemma}[theorem]{Lemma}
\newtheorem{corollary}[theorem]{Corollary}
\newtheorem{proposition}[theorem]{Proposition}

\newtheorem*{theoremA}{Theorem A}
\newtheorem*{theoremB}{Theorem B}
\newtheorem*{theoremC}{Theorem C}

\theoremstyle{definition}

\numberwithin{equation}{section}
\newcommand{\N}{{\mathbb{N}}}

\newcommand{\R}{{\mathbb{R}}}
\newcommand{\C}{{\mathbb{C}}}

\def\1{{{\mathit 1} \!\!\>\!\! I} }
\renewcommand{\phi}{\varphi}

\title[Global asymptotic stability]{Vector fields whose linearisation  \linebreak  is Hurwitz almost everywhere}

\subjclass[2000]{Primary 34D23; 37B25 Secondary 37C10}
\email{benito@ffclrp.usp.br}
\email{rrabanal@pucp.edu.pe}
\dedicatory{Dedicated in memoriam to professor Carlos Guti\'errez}

\medskip

\begin{document}

\maketitle

\centerline{\scshape Benito Pires \footnote{The first author was partially supported by FAPESP-BRAZIL (2009/02380-0 and \mbox{2008/02841-4)}}}
\medskip
{\footnotesize
 \centerline{Departamento de Computa\c c\~ao e Matem\'atica, Faculdade de Filosofia, Ci\^encias e Letras}
   \centerline{Universidade de S\~ao Paulo, Ribeir\~ao Preto, SP, Brazil}
} %

\medskip

\centerline{\scshape Roland Rabanal \footnote{The second author was partially supported by 
PUCP-PERU (DAI 2010-0058) and ICTP-ITALY (220 (Maths) RR/ab).}
}

\medskip

{\footnotesize
 \centerline{Secci\'on Matem\'atica, Pontificia Universidad Cat\'olica del Per\'u}
   \centerline{Av. Universitaria 1801, San Miguel: Apartado Postal 1761, Lima 100, Per\'u}
} %

\bigskip


\begin{abstract} A real matrix is Hurwitz if its eigenvalues have negative real parts. The following generalisation of the Bidimensional Global Asymptotic Stability Problem (BGAS) is provided:
 Let $X:\R^2\to\R^2$ be a $C^1$ vector field whose derivative $DX(p)$ is Hurwitz for almost all $p\in\R^2$. Then the singularity set of $X$, {\rm Sing}\,(X), is either an emptyset, a one--point set or a non-discrete set. Moreover, if ${\rm Sing}\,(X)$ contains a hyperbolic singularity then $X$ is topologically equivalent to the radial vector field $(x,y)\mapsto (-x,-y)$. This generalises BGAS to the case in which the vector field is not necessarily a local diffeomorphism.
  \end{abstract}

\maketitle

\section{Introduction}

This paper provides sufficient conditions for a planar $C^1$ vector field to be globally asymptotically stable, that is, to be topogically equivalent  to the radial vector field $(x,y)\mapsto (-x,-y)$. The strongest result in this respect is due to Guti{\'e}rrez \cite{G1}, Fessler \cite{F} and Glutsyuk \cite{Glu}
who solved positively the Markus-Yamabe Problem \cite{MY} (also known as the Bidimensional Global Asymptotic Stability Problem): 

\begin{theorem}[BGAS]\label{FG} Let $X:\R^2\to\R^2$ be a $C^1$ vector field $($map$)$ such that \mbox{$X(0)=0$}.
If the derivative $DX(p)$ is Hurwitz for all $p\in\R^2$ then $X$ is globally asymptotically stable $($globally injective$)$.
\end{theorem}

The condition that the derivative $DX(p)$ is Hurwitz is equivalent to requiring ${\rm Trace}\,(DX(p))<0$ and ${\rm Det}\,(DX(p))>0$. Therefore,
Theorem \ref{FG} can only be applied to vector fields that are local diffeomorphisms. 

It is not difficult to find globally asymptotically stable vector fields whose derivative is Hurwitz everywhere  except along curves or points where either the determinant or the trace of the derivative vanishes. For example, the derivative of the planar vector field
$X:\R^2\to\R^2$ defined by $X(x,y)=\left(-(x+1)^3+1,-(x+1)^2(y+1)+1\right)$ is Hurwitz Lebesgue almost every\-where but not everywhere. Besides, $X(0)=0$. We cannot apply \mbox{Theorem 1.1}, but we can apply \mbox{Theorem A }, to conclude that $X$ is globally asymptotically stable. Notice that this vector field (as a map) is neither a local diffeomorphism nor an injective map.
The aim of this paper is to provide a generalisation of  \mbox{Theorem \ref{FG}} which includes vector fields which are not necessarily local diffeomorphisms.

This paper also contains a result (Theorem B) about the cardinality of the singularity set of a $C^1$ vector field whose derivative is Hurwitz at almost all points of $\R^2$. It states that under such hypothesis the singularity set of the vector field has either zero, one or infinitely many points. 

Finally, we provide some applications of our results to gradient and Hamiltonian $C^1$ vector fields. In particular, we also give a sufficient condition (Theorem C) for a vector field to be a global center. This result turns out to be a generalisation of the paper of \cite{R} to the case in which $X$ is not necessarily a local diffeomorphism.\\

\noindent {\bf Acknowledgments.} This paper was developed during a visit of the authors to the Abdus Salam International Centre for Theoretical Physics (ICTP). We would like to thank the members of the Mathematics Section for their kind hospitaly. 
\section{Statement of the results}

Throughout this paper, ``for almost all" or ``a.e." means ``for Lebesgue almost all". We say that $p\in\R^2$ is a simple singularity of a $C^1$ vector field $X:\R^2\to\R^2$ if $X(p)=0$ and ${\rm Det}\,(DX(p))\neq 0$. We say that $p\in\R^2$ is a hyperbolic singularity of $X$ if $X(p)=0$ and the eigenvalues of $DX(p)$
have non-null real parts. By ``topologically equivalent" we understand ``globally topologically equivalent". 
For the sake of  simplicity we assume, whenever the singularity set of the vector field $X$ is not empty, that $X(0)=0$. 

The main results of this paper are the following:

\begin{theoremA} Let $X:\R^2\to\R^2$ be a $C^1$ vector field whose derivative $DX(p)$ is Hurwitz for almost all $p\in\R^2$. If $0\in\R^2$ is a hyperbolic singularity of $X$ then $X$ is topologically equivalent to the radial vector field $(x,y)\mapsto (-x,-y)$.
\end{theoremA}

\begin{theoremB} Let $X:\R^2\to\R^2$ be a $C^1$ vector field. If the derivative $DX(p)$  is Hurwitz for almost all
$p\in\R^2$ then the singularity set of $X$ is either an emptyset, a one-point set or a non-discrete set.

 \end{theoremB} 
 
 We say that a complex number $z\in\C$ is purely imaginary if $z=b i$ for some real number $b\neq 0$. Given $\epsilon\in \{-1,1\}$, we call the vector field $(x,y)\mapsto \epsilon (-y,x)$ a linear center.

 \begin{theoremC} Let $X:\R^2\to\R^2$ be a $C^1$ vector field whose derivative $DX(p)$ has purely imaginary eigenvalues for almost all $p\in\R^2$. If $0\in\R^2$ is a simple singularity of $X$ then $X$ is topologically equivalent to  a linear center.
  \end{theoremC}
  
  We also include the following interesting examples (see the proof in Section 6):

\begin{proposition}\label{example} The following holds true: 
\begin{itemize} 
\item [(a)] There exists a non-injective and globally asymptotically stable \mbox{$C^1$} vector field whose derivative $DX(p)$ is Hurwitz for almost all $p\in\R^2$. 
\item [(b)] There exists a $C^1$ vector field whose derivative $DX(p)$ is Hurwitz for almost all $p\in\R^2$ and whose singularity set is a line $($a non-discrete set$)$.
\end{itemize}
\end{proposition}

To explain these results, let us to fix the following notation. We say that a $C^1$ vector field $X:\R^2\to\R^2$ is Hurwitz (respectively Hurwitz a.e.) if the derivative $DX$ is Hurwitz everywhere (respectively Hurwitz almost everywhere).

If the hypothesis of Hurwitz of Theorem \ref{FG} is replaced by Hurwitz a.e. one cannot ensure anymore
that $0$ is a global attractor of $X$. In fact, we present an example of a vector field Hurwitz a.e. which has the line \mbox{$\{(x,y)\in\R^2\mid x=-1\}$} as its singularity set. However, if we rule out exotic behaviour by requiring that $0$ is a hyperbolic singularity then the hypothesis of Hurwitz a.e. still implies that $0$ is a global attractor (Theorem A).  

Theorem B characterizes the singularity set of a Hurwitz a.e. vector field: it has either zero, or  one or infinitely many points. In particular, a vector field having exactly two singularities cannot be Hurwitz a.e. 

Theorem C gives spectral conditions for a $C^1$ vector field to be a global center. It was proved in
\cite{R} that if a differentiable vector field $X:\R^2\to\R^2$ has the property that $X(0)=0$ and $DX(p)$ has purely imaginary eigenvalues for all $p\in\R^2$ then $0$ is a global center of $X$. Theorem C says that this result is still true for
a vector field whose derivative $DX(p)$ has purely imaginary eigenvalues for almost all $p\in\R^2$, provided $0\in\R^2$ is a simple singularity of $X$.

It is also included in this paper some applications (Corollaries \ref{hvf} and \ref{gvf}) of the main theorems to gradient and Hamiltonian $C^1$ vector fields. A consequence of these corollaries is that the Jacobian determinant of a gradient (or Hamiltonian) $C^1$ vector field that has more than one singularity changes signal. This is exploited in the section of examples.  

In addition to the articles mentioned so far, we would like to attach some important historical remarks on the theme of global asymptotic stability and injectivity. Olech-Meisters \cite{MOMO} gave an affirmative answer to the Markus-Yamabe Problem (in the plane) in the polynomial case. In \cite{O}, Olech  pointed out the relations between global asymptotic stability and global injectivity. The problem remained opened for more than 30 years until it was completely solved by Gutierrez \cite{G1}, Fessler \cite{F} and Glutsyuk \cite{Glu}  in an affirmative way for $C^1$ vector fields. A generalisation of this result for differentiable (not necessarily $C^1$) vector fields was given in
 Fernandes-Gutierrez-Rabanal \cite{FGR2} (see also \cite{GPR}). Concerning the dynamics of planar Hamiltonian vector fields, we refer the reader to the article of Jarque-Nitecki \cite{JN}.  

\section{Vector fields with Hurwitz derivatives a.e.}

Given a $C^1$ vector field $X:\R^2\to\R^2$, we let:
$${\rm Spc}{\,(X)}=\{\text{Eigenvalues of}\,\, DX(p)\mid p\in \R^2\}\subset\C.$$

\begin{proposition}\label{U} Let $X:\R^2\to\R^2$ be a $C^1$ vector field such that $DX(p)$ is  Hurwitz for almost all $p\in\R^2$. Let $U=\{p\in\R^2\mid {\rm Det}\,(DX(p))>0\}$. Then $U$ is an open, dense subset of $\R^2$ and $X\vert_U$ is injective.
 \end{proposition}
\begin{proof} The hypothesis that $DX$
is Hurwitz almost everywhere yields the  inequalities \mbox{${\rm Trace}\,(DX(p))<0$} and ${\rm Det}\,(DX(p))>0$
for almost every $p\in\R^2$. Now by the continuity of the partial derivatives of $X$, we have that $U$ is open and dense, 
\mbox{${\rm Trace}\,(DX(p))\le 0$} and ${\rm Det}\,(DX(p))\ge 0$ for all $p\in\R^2$. 
Thus if we set $\epsilon_n=\frac1n$ and $X_n=X-\epsilon_n I$ then $${\rm Spc}\,(X_n)\subset \{z\in\C\mid {\rm Re}\,(z)<0\}$$
 for all $n\in\N$, where ${\rm Re}\,(z)$ stands for the real part of the number $z$.
  By \mbox{Theorem \ref{FG}}, we have that $X_n$ is injective for all $n\in \N$. We claim that $X\vert_U$ is injective. Choose $p,q\in U$ such that $X(p)=y=X(q)$. We will prove that $p=q$. By the definition of $U$ and by the Inverse Function Theorem, there exist compact neighborhoods, $U_1,U_2,V$ of $p,q,y$, respectively, such that
$X\vert_{U_i}:U_i\to V$ is a homeomorphism and $U_1\cap U_2=\emptyset$. By the definition of $X_n$,
for $n$ large enough, $X_n(U_1)\cap X_n(U_2)$ will contain a neighborhood $W$ of $y$. Hence,
for all $w\in W$, $\#(X_n^{-1}(w))\ge 2$. This contradicts the injectivity of $X_n$.
\end{proof}

\begin{corollary}\label{2e} Let $X:\R^2\to\R^2$ be a $C^1$ vector field such that $DX(p)$ is  Hurwitz for almost all $p\in\R^2$. If $0\in\R^2$ is a hyperbolic singularity of $X$ then for all $\rho>0$ there exists $\epsilon>0$ such that $\Vert X(p)\Vert> \epsilon$ for all $\Vert p\Vert > \rho$. In particular, ${\rm Sing}\,(X)=\{0\}$.
\end{corollary}
\begin{proof} It suffices to prove the result for arbitrarily small $\rho>0$.
Let $U$ be as in Proposition \ref{U}. 
Then $0\in U$ and $X(0)=0$. Let $\rho>0$ be such that $B_{\rho}\subset U$, where $B_\rho$ an open ball of ratio $\rho$ centered at $0$. Because $X\vert_U$
is an open map and $X(0)=0$, there exists a ball $B_{2\epsilon}$ of ratio $2 \epsilon$ centered at $0$ such that   $0\in B_{2 \epsilon}\subset X(B_\rho)$. It follows from the injectivity of $X\vert_U$ that $X(U\setminus B_\rho)\cap B_{2\epsilon}=\emptyset$. Thus, $\Vert X(p)\Vert>2 \epsilon$ for all $p\in U\setminus B_\rho$. We claim that $\Vert X(p)\Vert> \epsilon$ for all $p$ such that $\Vert p\Vert>\rho$. In fact, if $p\in\R^2\setminus B_\rho$ then, by the density of $U$, there exists a sequence $\{p_n\}_{n=0}^\infty\subset U\setminus B_\rho$ tending to $p$ such that $\Vert X(p)\Vert=\lim_{n\to\infty} \Vert X(p_n)\Vert$. As $p_n\in U\setminus B_\rho$, we have $\Vert X(p_n)\Vert> 2 \epsilon$ for all $n$ and hence $\Vert X(p)\Vert\ge 2\epsilon>\epsilon$.
\end{proof}

Given a $C^1$ vector field $X=(f,g):\R^2\to\R^2$, let $X^*=(-g,f)$ be the orthogonal
vector field to $X$. The same notation as that for intervals of $\R$
will be used for oriented arcs of trajectory $[p,q], [p,q)$,...
(respectively $[p,q]^*,[p,q)^*$,...) of $X$ (respectively $X^*$), connecting the
points $p$ and $q$. The orientation of theses arcs is that induced
by $X$ (respectively $X^*$).  

 A {\it compact rectangle} $R=R(p_1,p_2;q_1,q_2)\subset \R^2$ of a $C^1$ vector field $X:\R^2\to\R^2$
 is the compact region the boundary of which is made up of two arcs of
trajectory $[p_1,p_2], [q_1,q_2]$ of $X$ and two arcs of trajectory
$[p_1,q_1]^*,[p_2,q_2]^*$ of $X^*$. Notice that we assume that the
flow induced by $X$ goes into $R$ by $[p_1,q_1]^*$ and leaves $R$ by
$[p_2,q_2]^*$.

For any arc of trajectory $[p,q]^*$ of $X^*$, let
$$L([p,q]^*)=\Big\vert\int_{[p,q]^*}\Vert X^*\Vert\,ds\Big\vert,$$
where $ds$ denotes the arc length element. Given an arc of
trajectory [p,q] (resp. $[p,q]^*$), we denote by $\ell([p,q])$
(resp. $\ell([p,q]^*)$) the arc length of it. 

\begin{lemma}{\label{Green}} Let $R=R(p_1,p_2;q_1,q_2)\subset \R^2$ be a compact rectangle of a $C^1$ vector field. Then
$$L([p_2,q_2]^*)-L([p_1,q_1]^*)=\int_R {\rm Trace}\,{(DX)}\,dx\wedge dy. $$
\end{lemma}
\begin{proof} See Green's Formula as presented in \cite[Corollary 5.7]{pfe}.
\end{proof}

 In what follows, we let $\gamma_p^+$ denote a maximal solution of the ordinary differential equation
$u'=X(u)$ passing through $p\in\R^2$. We let $[0,b_+)$, $b_+\in [0,+\infty)\cup\{+\infty\}$, denote the maximal interval of definition of $\gamma_p^+$. We identify $\gamma_p^+$ with its image $\gamma_p^+\left([0,b_+)\right)$ and we call
$\gamma_p^+$ a positive semitrajectory of $X$. Given an open set $V$, we set 
$d(\gamma_p^+,V)=\inf \{\Vert \gamma_p^+(t)-v\Vert: t\in [0,b_+),v\in V  \}$. 
We denote by $\omega(p)$ the $\omega$-limit set of a point $p$:
$$
\omega(p)=\left\{q\in\R^2\mid \exists t_n\to b_+\,\,\text{such that}\,\, \lim_{n\to\infty} \gamma_p^+(t_n)=q        \right\}.
$$

An alternative proof of next proposition may be found in \cite[Lemma 3.4, p. 480]{FGR2}. A complete proof is provided here.

\begin{proposition}\label{ei} Let $X:\R^2\to
\R^2$ be a $C^1$ vector field whose derivative $DX(p)$ is Hurwitz for almost all $p\in\R^2$. If $0\in\R^2$ is a hyperbolic singularity of $X$
 then the set \mbox{$\{p\in\R^2\mid \omega(p)=\emptyset\}$} is open.
\end{proposition}
\proof 
Let $p_1\in\R^2$ be such that $\omega(p_1)=\emptyset$. 
We will prove that $\omega(p)=\emptyset$ for all $p$ in a neighborhood of $p_1$. Let $B_r=\{p\in\R^2: \Vert p\Vert <r\}$. 
Because $\omega(p_1)=\emptyset$, there exists $K>0$ such that $\gamma_{p_1}^+\subset \R^2\setminus B_{2K}$.
Let $V=\{p\in\R^2\mid d(p,\gamma_{p_1}^+)<K/4\}$ be an open neighborhood of $\gamma_{p_1}^+$ of ratio $K/4$. Let $W\subset V\subset\R^2\setminus B_{2K}$ be the union of all arcs of trajectory of 
 $X^*$ contained in $V$ and intersecting $\gamma_{p_1}^+\setminus\{p_1\}$. 
 By the Long Tubular Flow Theorem \cite[Proposition 1.1, p. 93]{PM},
 $W$ is an open neighborhood of $\gamma_{p_1}^+\setminus\{p_1\}$. By Corollary \ref{2e}, there exists $\epsilon>0$ such that
 such that $\Vert X(p)\Vert>\epsilon$ for all $p\in \R^2\setminus B_{K}$. Let $[x_1,y_1]^*\subset \overline{W}\subset V$ be an arc of trajectory of
 $X^*$ containing $p_1$ such small that  $\frac{d_1}{\epsilon}\ell([x_1,y_1]^*)<\frac{K}{8}$, where $d_1=\sup\,\{\Vert X(w)\Vert:w\in [x_1,y_1]^*\}$.
 We claim that $\gamma_{q_1}^+\subset V$ for all $q_1\in [x_1,y_1]^*$. Suppose that the Claim is false. 
So, without loss of generality, we may assume that $\gamma_{{q_1}}^+\not\subset V$ for some $q_1\in [p_1,y_1]^*$. Then there exists $\eta>0$ such that $\gamma_{q_1}^+(t)\in W$ for all $t\in (0,\eta)$ but $q_2:=\gamma_{q_1}^+(\eta)\not\in W$.
 In  this way, $(q_1,q_2)\subset W$. Hence, for all $q\in (q_1,q_2)$, there exists $p\in \gamma_{p_1}^+$ such that
 $[p,q]^*$ is an arc of trajectory of $X^*$.
 By Lemma \mbox{\ref{Green}}, for the rectangle
$R:=R(p_1,p;q_1,q)$ we have:
$$L([p,q]^*)-L([p_1,q_1]^*)=\int_R{\rm Trace}\,(DX)\,dx\wedge dy<0.
$$
As $R\subset V$ and $\Vert X(p)\Vert>\epsilon$ for all $p\in V$, we have that:
\begin{eqnarray*}
\epsilon\ell([p,q]^*)\le\Big\vert\int_{[p,q]^*}\Vert X\Vert ds\Big\vert &=&
L([p,q]^*)<L([p_1,q_1]^*)=\\&=& \Big\vert\int_{[p_1,q_1]^*}\Vert
X\Vert ds\Big\vert\le d_1\ell([p_1,q_1]^*).
\end{eqnarray*}
Therefore:
$$\ell([p,q]^*)\le\frac{d_1}{\epsilon}\ell([p_1,q_1]^*)\le\frac{d_1}{\epsilon}\ell([x_1,y_1]^*)<\frac{K}{8}.
$$
We have shown above that for each $t\in (0,\eta)$ there exist $q(t)=\gamma_{q_1}^+(t)\in (q_1,q_2)$ and $p(t)\in \gamma_{p_1}^+$ such that $[p(t),q(t)]^*$ is an arc of trajectory of $X^*$ of arc length less than $K/8$. It is easy to see that there exists $p_2\in \gamma_{p_1}^+$ such that $\lim_{t\to\eta} p(t)=p_2$.
In this way, the limit points of the arc of trajectory $[p(t),q(t)]^*$ of $X^*$ is an arc of trajectory of $X^*$
that connects $p_2$ and $q_2$, denoted by $[p_2,q_2]^*$.
This means that either $q_2\in W$ (impossible by the choice of $q_2$)
or $q_2\in \partial V$ (impossible because $\ell([p_2,q_2]^*)<K/8<K/4$). Thus 
 $\gamma_{q_1}^+\subset V$ for all $q_1\in [x_1,y_1]^*$. Consequently, $\omega(q_1)=\emptyset$ for all
 $q_1\in [x_1,y_1]^*$.
 By the Long Tubular Flow Theorem, if we take a small ball $B(p)\subset V$ centered at $p$, all the trajectories $\gamma_p^+$, with $p\in B(p)$ will cross the transverse segment
$[x_1,y_1]^*$. Thus $\omega(p)=\emptyset$ for all $p\in B(p)$.
\endproof

\begin{theoremA} Let $X:\R^2\to\R^2$ be a $C^1$ vector field whose derivative $DX(p)$ is Hurwitz for almost all $p\in\R^2$. If $0\in\R^2$ is a hyperbolic singularity of $X$ then $X$ is topologically equivalent to the radial vector field $(x,y)\mapsto (-x,-y)$.
\end{theoremA}
\begin{proof}
 By Corollary \ref{2e}, $0\in\R^2$ is the unique singularity of $X$. Now let
$$W^s=\{p\in\R^2\mid \omega(p)=\{0\}\}.$$
It is easy to show that the eigenvalues of $DX(0)$ have negative real parts.
Thus $0$ is a local attractor so that $W^s$ is not empty. Besides, by the Long Flow Box Theorem, $W^s$ is open. It follows from Green's Formula that $X$ has no periodic trajectories. The Poincar\'e-Bendixson Theorem implies that for each $p\in\R^2$,
either $\omega(p)=\{0\}$ or $\omega(p)$ is an empty set. According to Proposition \ref{ei}, we have that
$\R^2\setminus W^s$ is an open set. In this way, $W^s$ is an open and closed subset of $\R^2$ and so $W^s=\R^2$.
\end{proof}

\section{Singularity set of Hurwitz a.e. vector fields}

\begin{theoremB} Let $X:\R^2\to\R^2$ be a $C^1$ vector field. If the derivative $DX(p)$  is Hurwitz for almost all
$p\in\R^2$ then the singularity set of $X$ is either an emptyset, a one-point set or a non-discrete set.
 \end{theoremB}
 \begin{proof} 
 
 Let $S={\rm Sing}\,(X)$ be the singularity set of $X$. Suppose that $S$ is discrete. We have to show that $S$ is either empty or a one-point set. Thus we may assume that $S$ is a non-empty discrete set.
 Consequently, we may write $S=\{z_i\}$, where $X(z_i)=0$ for all $i$. The claims below show that under the above hypotheses, $S$ is a one-point set.\\
 
 \noindent Claim 1: Let $z\in S$ be fixed. There exist a sequence $\{X_k\}_{k=1}^\infty$ of $C^1$ vector fields
  tending to $X$ in
 the $C^1$ uniform topology and a sequence of points $\{v_k\}_{k=1}^\infty\subset\R^2\setminus S$ tending to $z$ such
 that $X_k(v_k)=0$ and $DX_k(v_k)$ is Hurwitz for all $k$. Furthermore, each $X_k$ is Hurwitz a.e., that is, its derivative $DX_k(p)$ is Hurwitz for almost all $p\in\R^2$.\\
 
 \noindent {\it Proof of Claim 1}. Without loss of generality, we may assume that $z=0$. Let $B_k=\{p\in\R^2:\Vert p\Vert<\frac1k\}$. For each $k\in\N$, by the continuity of $X$ and because $X(0)=0$, there exists $n_k>k$ such that $X(B_{n_k})\subset B_k$. For each $k\in\N$, let $v_k\in B_{n_k}\setminus S$ be such that $DX(v_k)$
 is Hurwitz. Define $\omega_k=X(v_k)$. Notice that $\Vert \omega_k\Vert=\Vert X(v_k) \Vert<\frac1k$ and $\Vert v_k\Vert<\frac{1}{n_k}<\frac1k$. For each $k\in\N$, set $X_k=X-\omega_k$. It is plain that $X_k(v_k)=0$ and that $DX_k(v_k)=DX(v_k)$ is Hurwitz. In particular, $v_k$ is a simple singularity of $X_k$. As $X(p)-X_k(p)=w_k$ for all $p\in\R^2$, we have that $X_k$ tends to $X$ in the $C^1$ uniform topology. Finally, we have that $DX_k(p)=DX(p)$ for all $p\in\R^2$. Hence, $DX_k(p)$ is Hurwitz for almost all $p\in\R^2$ and for each $k\in\N$. 
 \\
 
\noindent Claim 2: Each singularity $z\in S$ has index 1 .\\

\noindent {\it Proof of Claim 2}. Let $z\in S$ be fixed. Let $\{X_k\}_{k=1}^\infty$ and $\{v_k\}_{k=1}^\infty$ be as in \mbox{Claim 1}. Given a topological circle $\gamma$ contained in $\R^2$, we call  the open bounded connected component
of $\R^2\setminus \gamma$ the interior of $\gamma$, denoted by ${\rm Int}\,(\gamma)$.
Now let $\gamma\subset\R^2$ be a topological circle such that $z\in {\rm Int}\,(\gamma)$ and ${\rm Int}\,(\gamma)\cap \left(S\setminus\{z\}\right)=\emptyset$.
As $v_k\to z$ as $k\to \infty$, we may assume that $v_k\in {\rm Int}\,(\gamma)$ for all $k$. 
By Theorem A, we have that
each $X_k$ is globally asymptotically stable. In this way, the index of $X_k$ along $\gamma$ is equal to 1. By continuity, as $X_k\to X$ uniformly, we have that the index of $X$ along $\gamma$ is 1. As $z$ is the only singularity of $X$ contained in ${\rm Int}\,(\gamma)$ we have that its index is 1.  \\

\noindent Claim 3: $S$ is a one-point set. \\
 
 \noindent {\it Proof of Claim 3}. Suppose that $S$ has at least two singularities, say $z_1$ and $z_2$.
Let $\gamma$ be a topological circle
 such that $\{z_1, z_2\}\subset {\rm Int}\,(\gamma)$ and ${\rm Int}\,(\gamma)\cap(S\setminus\{z_1,z_2\} )=\emptyset$. We can use the reasoning in the 
 proof of  Claim 2 to show that the index of the vector field along $\gamma$ is 1. This contradicts Claim 2, which states that ${\rm Int}\,(\gamma)$ contains exactly two singularities of index 1. So $S$ is a one-point set and the proof is finished.
 \end{proof}
 
 \section{Vector Fields with non-hyperbolic derivatives a.e.}
 
 In this section we prove Theorem C. We will need some lemmas.
 
 \begin{proposition}\label{tcca} Let $X:\R^2\to\R^2$ be a $C^1$ vector field whose derivative $DX(p)$ has purely imaginary eigenvalues for almost all $p\in\R^2$. Let $U=\{p\in\R^2\mid {\rm Det}\,(DX(p))>0\}$. Then $U$ is an open, dense subset of $\R^2$ and $X\vert_U$ is injective. 
 \end{proposition}
 \begin{proof} The same proof of Proposition \ref{U} holds word by word. The only change is to replace
 ${\rm Trace}\,(DX(p))<0$ by ${\rm Trace}\,(DX(p))=0$ a.e.
 \end{proof}
 
 \begin{corollary}\label{44} Let $X:\R^2\to\R^2$ be a $C^1$ vector field whose derivative $DX(p)$ has purely imaginary eigenvalues for almost all $p\in\R^2$. If $0\in\R^2$ is a simple singularity of $X$ then for all $\rho>0$ there exists $\epsilon>0$ such that $\Vert X(p)\Vert> \epsilon$ for all $\Vert p\Vert > \rho$. In particular, 
${\rm Sing}\,(X)=\{0\}$.

\end{corollary}
\begin{proof} It follows from Proposition \ref{tcca}. The proof is as in Corollary \ref{2e}.  
 \end{proof}
  
 \begin{theoremC} Let $X:\R^2\to\R^2$ be a $C^1$ vector field whose derivative $DX(p)$ has purely imaginary eigenvalues for almost all $p\in\R^2$. If $0\in\R^2$ is a simple singularity of $X$ then $X$ is topologically equivalent to a linear center.
  \end{theoremC}
  \begin{proof} We can adapt the results presented in \cite{R} to prove Theorem C. We point out here just the key steps. By Corollary \ref{44} and by the proof of \cite[Proposition 2.3, p. 656 ]{R}, $X$ has no hyperbolic sectors at infinity. Because $0$ is a simple singularity of $X$ with positive determinant, we have that the singularity $0$ of $X$ has index 1. Given a topological circle $\gamma$, let ${\rm Int}\,(\gamma)$ denote the bounded connected component of $\R^2\setminus\gamma$. By the above, the index of $X$ along any topological circle $\gamma$ such that $0\in {\rm Int}\,(\gamma)$ is 1. We claim that there exists a sequence of periodic trajectories $\{\Gamma_i\}_{i=1}^\infty$ such that  ${\rm Int}\,(\Gamma_i)\subset {\rm Int}\,(\Gamma_{i+1})$ for all $i\ge 1$ and $\cap_{i=1}^\infty \,\left(\R^2\setminus{\rm Int}\,(\Gamma_i)\right)=\emptyset$. Suppose that the claim is false. Let $\gamma$ be a topological circle, piecewise $C^1$ smooth, which has the least number of topological tangencies with $X$ and such that $0\in{\rm Int}\,(\gamma)$. By the above, as {\rm Sing}\,(X)=\{0\}, the index of $X$ along $\gamma$ is 1. Because there is no hyperbolic sector at infinity nor singularities in $\R^2\setminus {\rm Int}\,(\gamma)$, it is possible to construct a transverse circle $C$ to $X$ such that $\gamma\subset {\rm Int}\,(C)$. This contradicts Green's Theorem if we use the fact that \mbox{${\rm Trace}\,(DX)=0$ a.e}. Hence, the claim is true.
  Around the origin, the fact that $0$ is a simple singularity such that $DX(0)$ has purely imaginary eigenvalues implies that $0$ is a local center of $X$.
  This together with the claim and the hypothesis that ${\rm Trace}\,(DX)=0$ a.e. leads to the conclusion, after applying some standard arguments of topological dynamics, that $X$ is topologically equivalent to a linear center.
   \end{proof}
 
 \section{Corollaries of the main theorems}
 
 We say that $X:\R^2\to\R$ is a $C^1$ Hamiltonian vector field if there exists a \mbox{$C^2$ function}
 $h:\R^2\to\R$ such that $X=(-h_y,h_x)$. 
 
 \begin{corollary}\label{hvf} Let $X:\R^2\to\R^2$ be a Hamiltonian $C^1$ vector field whose Jacobian determinant is positive almost everywhere. If  $0\in\R^2$ is a simple singularity of $X$ then $X$ is topologically equivalent to a linear center.
 \end{corollary} 
 \begin{proof} 
By the Schwarz Lemma, the mixed partial derivatives of $X$ are equal to each other, $h_{yx}=h_{xy}$. Hence, ${\rm Trace}\,(DX(p))=-h_{yx}(p)+h_{xy}(p)=0$ for all $p\in\R^2$. By the hypothesis on the Jacobian of $X$, we have that ${\rm Det}\,(DX(p))>0$ for almost all $p\in\R^2$. Putting it all together, $DX(p)$ has purely imaginary eigenvalues for almost all $p\in\R^2$. Now the result follows from Theorem C.
\end{proof}
 
 We say that a vector field $X:\R^2\to\R^2$ is $C^1$ gradient if there exists a \mbox{$C^2$ function} $h:\R^2\to\R$ such that $X=(h_x,h_y)$. We say that a vector field $Y:\R^2\to\R^2$ is dissipative if ${\rm Trace}\,(DY(p))<0$ for almost all $p\in\R^2$. If $X$ is a $C^1$ gradient vector field this does not mean that $X$ or $-X$ is dissipative (for example, let $X:\R^2\to\R^2$ be defined by $X=(3x^2,-3y^2)$).
 
 \begin{corollary}\label{gvf} Let $X:\R^2\to\R^2$ be a gradient $C^1$ vector field whose Jacobian determinant is positive almost everywhere. If $0\in\R^2$ is a hyperbolic singularity of $X$  then $X$ is topologically equivalent to the radial vector field $(x,y)\mapsto (-x,-y)$.
 \end{corollary} 
 \begin{proof} Let $Y=RX$ be the vector field obtained from $X$ by the rotation $R:\R^2\to\R^2$ defined by $R(x,y)=(-y,x)$. It is plain that $X$ is a Hamiltonian vector field. By Corollary \ref{hvf}, $Y$ is a global center. This means that we can split $\R^2\setminus \{0\}$ as the union of transverse circles to $X$. It is not diffcult to show that $X$  is topologically equivalent to the radial vector field $(x,y)\mapsto (-x,-y)$.  
  \end{proof}

\section{Commented examples}

In this section we present a study of the spectrum (i.e. the distribution of the eigenvalues of $DX(p)$ as $p$ runs in $\R^2$) of many maps and vector fields. We group the cases with the same behaviour into subsections.

\subsection{Non-injective maps whose derivative is Hurwitz almost everywhere}

The Table 1 presents three examples of non-injective planar maps. The first map will be used to proof  the following proposition.\\

\noindent{\bf Proposition 2.1.}\hspace{-0.1cm} {\it The following holds true: 
\begin{itemize} 
\item [(a)] There exists a non-injective and globally asymptotically stable \mbox{$C^1$} vector field whose derivative $DX(p)$ is Hurwitz for almost all $p\in\R^2$. 
\item [(b)] There exists a $C^1$ vector field whose derivative $DX(p)$ is Hurwitz for almost all $p\in\R^2$ and whose singularity set is a line $($a non-discrete set$)$.
\end{itemize} } 
\begin{proof} (a) Let $X:\R^2\to\R^2$ be the $C^1$ vector field defined by $X(x,y)=F(x,y)=(-(x+1)^3+1,-(x+1)^2 (y+1) +1)$.
By the Table 1, we have that $DX(x,y)$ is Hurwitz everywhere except for the line $\{(x,y)\mid x=-1\}$. Furthermore, 
$p=(0,0)$ is a hyerbolic singularity of $X$. It follows from Theorem A that $X$ is globally asymptotically stable. This map is not injective because it takes the line $\{(x,y)\mid x=-1\}$ onto the point $(1,1)$. (b) The vector field $Y:\R^2\to\R^2$ defined by $Y(x,y)=(-x^3,-x^2y)$ has the line $S=\{(x,y)\mid x=0\}$ as its singularity set. Moreover, $DX(p)$ is Hurwitz for all $p\in\R^2\setminus S$.
\end{proof}

Olech \cite{O} proved that, for planar vector fields whose derivative is Hurwitz everywhere, injectivity of the vector field (considered as a map) and global asymptotically stability are equivalent. Proposition \ref{example} shows that such correspondence is lost under the weaker hypothesis that the derivative is Hurwitz almost everywhere.

\begin{table}[h]
\caption{Non-injective maps whose derivative is Hurwitz a.e.}
\scalebox{0.8}{%
\begin{tabular}{l | l | l }
\phantom{aaaaaaaaaaaaa} Map           & ${\rm Trace}\, \left(DF(x,y)\right)$ & $ {\rm Det}\,(DF(x,y))$  \\
\hline & & \\
$F(x,y)=(-(x+1)^3+1,-(x+1)^2 (y+1) +1)$& $-4(x+1)^2$  & $3(x+1)^4$ \\ 
& & \\
$G(x,y)=\left( \dfrac{-x^3}{1+x^2},\dfrac{-y x^2}{1+x^2}\right)$           & $\dfrac{-2 x^2 (2+x^2)}{(1+x^2)^2}$ &    $\dfrac{x^4 (3+x^2)}{(1+x^2)^3}$   \\ & & \\

$H(x,y)=\left[-1+\dfrac{2}{\pi} {\rm Tan}^{-1} \left( \dfrac{y}{x} \right) \right ]\cdot (x,y)$          & $-2+\dfrac{4}{\pi} {\rm Tan}^{-1} \left( \dfrac{y}{x} \right)  $& $\dfrac{\left(\pi-2{\rm Tan}^{-1}\left(\dfrac{y}{x} \right) \right)^2}{\pi^2}$ \\ & & \\
\hline
\end{tabular}}
\end{table}

As for the other examples in Table 1, $G$ is a rational map whose derivative is Hurwitz everywhere except for the line $\{(x,y)\mid x=0\}$. This set is also the set of zeros of $G$. Notice that by Proposition \ref{U}, $G$ is injective on $\R^2\setminus \{(x,y)\mid x=0\}$.   

The map $H$ can be continuously extended in a non-injective way to the whole plane. Its derivative is Hurwitz almost everywhere. Considered as vector fields, $G$ and $H$ are not globally asymptotically stable because their singularity sets contain a line of singularities. 

\begin{table}[h]
\caption{Hamiltonian Vector Fields}
\scalebox{0.85}{%
\begin{tabular}{l|l}
\phantom{aaaaaaaaaaaaaa}Vector Field            &    \phantom{aaaaaaaaa}  ${\rm Det}\,(DF(x,y))$  \\
\hline  
& \\
$X(x,y)=\left(-(y-1)^3-1,(x-1)^3+1 \right )$  & $ 9 (x-1)^2 (y-1)^2 $ \\ 
& \\ $Y(x,y)=\left(-2 e^{-(x^2+y^2)} x y, (2x^2-1) e^{-(x^2+y^2)}\right )$  & $-4 (2 x^4 + y^2 + x^2 (-3 + 2 y^2)) e^{-2 (x^2 + y^2)} $ \\ 
 & \\  $Z(x,y)= (-2 y + 4 y^3, -2 x + 4 x^3)$&  $ -4(-1+6 x^2)(-1+6y^2)$\\ & 
 \\  \hline
\end{tabular}}
\end{table}

\begin{figure}[h]
  \centering
  \subfloat[Vector field X]{\label{fig:gull}\includegraphics[width=4.3cm]{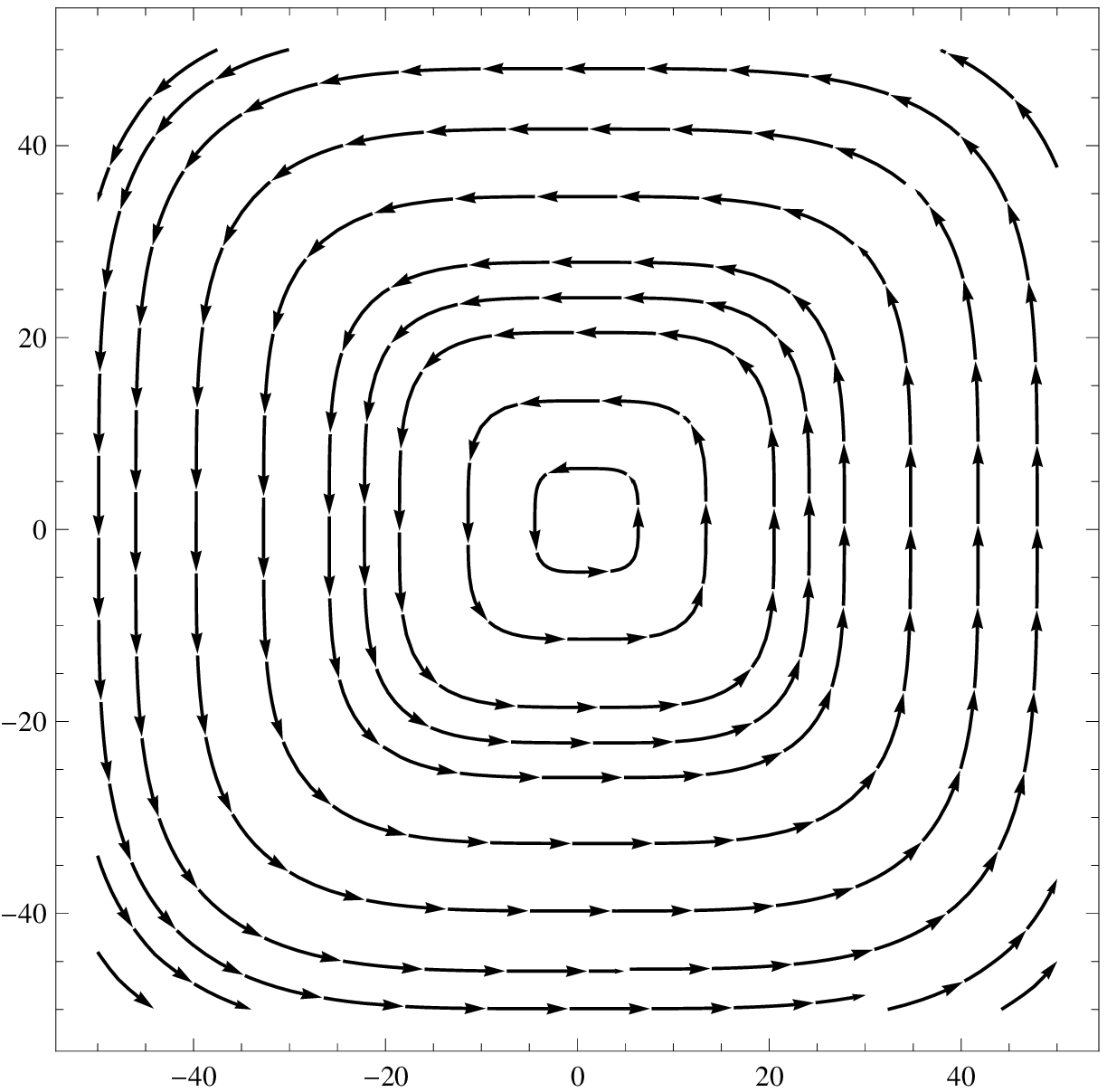}}                
  \subfloat[Vector field Y]{\label{fig:tiger}\includegraphics[width=4.3cm]{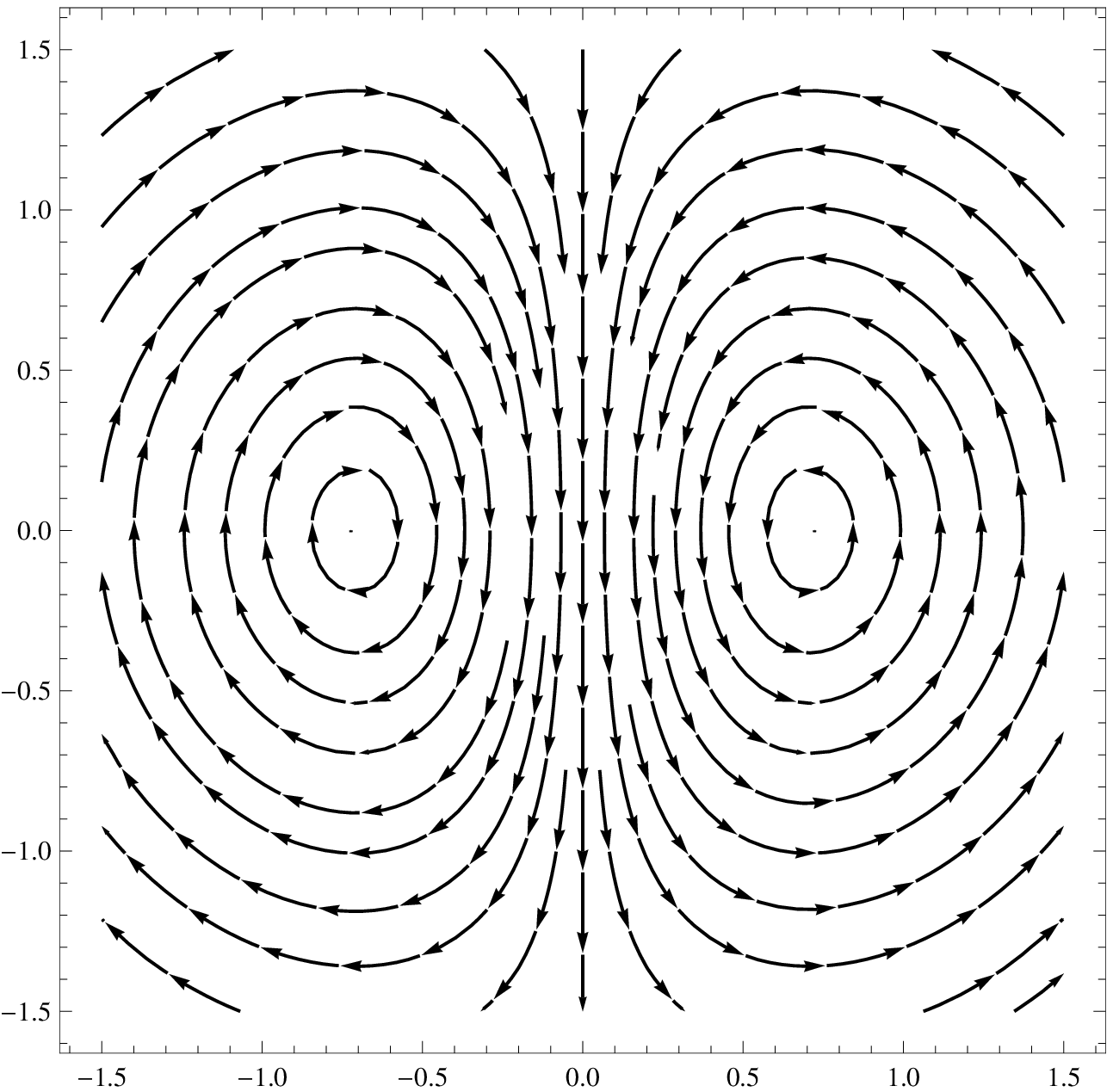}}
  \subfloat[Vector field Z]{\label{fig:mouse}\includegraphics[width=4.3cm]{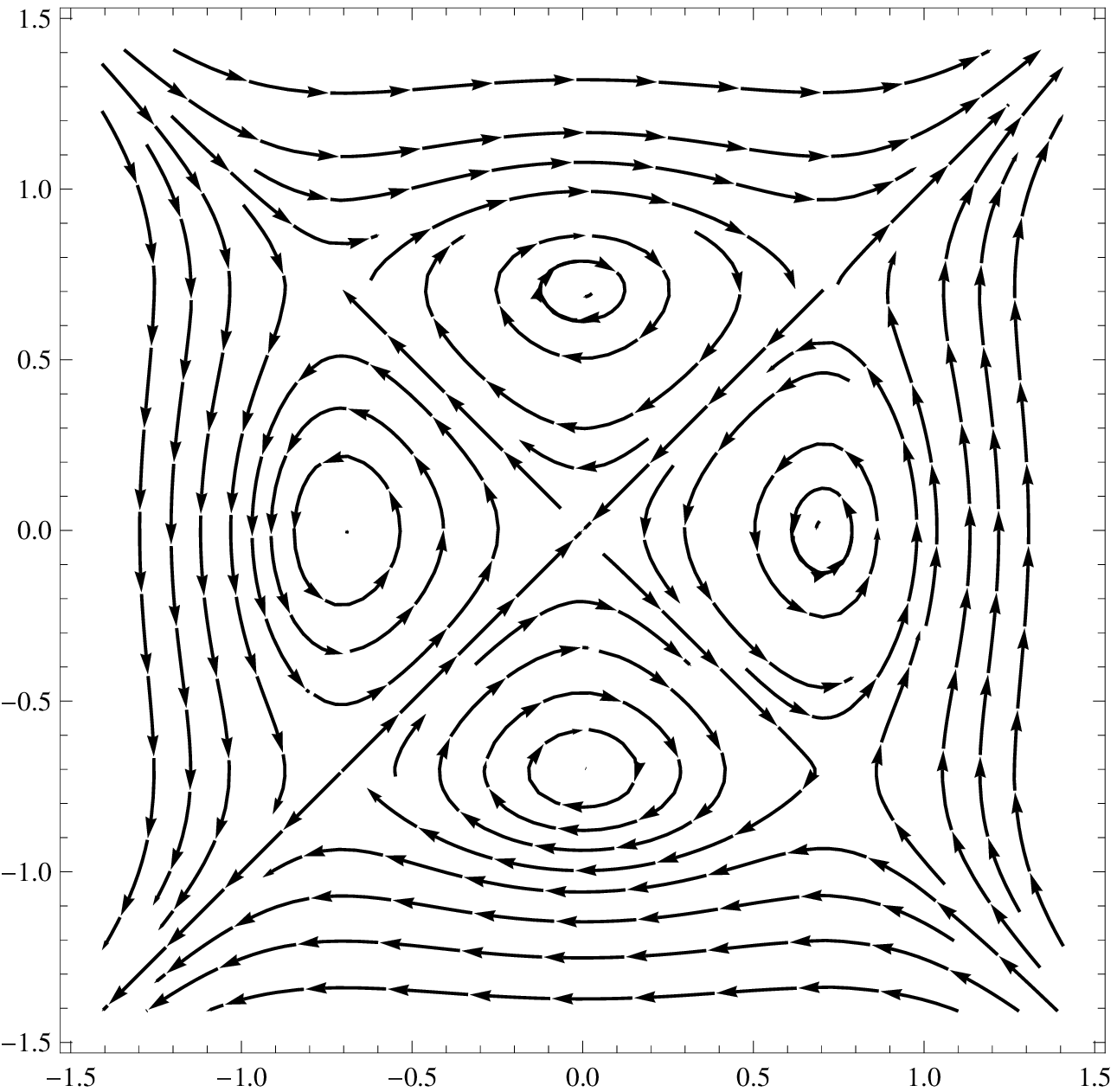}}
  \caption{Phaseportraits of Hamiltonian vector fields}
  \label{T2}
\end{figure}

\subsection{Hamiltonian vector fields and their Jacobians}

Corollary \ref{hvf}  may be applied to understand the phaseportrait of Hamiltonian $C^1$ vector fields. In \mbox{Table 2}  three examples of Hamiltonian vector fields are displayed. The first one, $X$,   
has the Jacobian positive almost everywhere. Besides, $0$ is a simple singularity of $X$. 
 In this way, By Corollary \ref{hvf}, $X$ is a global center  (see Figure \ref{T2}.A). Concerning the vector fields $Y$ and $Z$, as they have more than one singularity each, their Jacobian determinants
 have necessarily to change signal (see Figures \ref{T2}.B and \ref{T2}.C). This prescription agrees with the formula of the corresponding Jacobian determinants (Table 2).

\end{document}